\documentclass[12pt,reqno]{amsart}
\usepackage{latexsym,amsmath,amsthm,amssymb,amsfonts}
\usepackage{setspace}
\usepackage{color}
\usepackage{tikz}
\numberwithin{equation}{section}
\newtheorem{theorem}{Theorem}[section]
\newtheorem{lem}[theorem]{Lemma}

\newtheorem{pro}[theorem]{Proposition}

\newtheorem{defi}[theorem]{Definition}
\newtheorem{rem}[theorem]{Remark}

\def\S2{\mathbb{S}^2}
\def\s{\,\,\,\,}

\def\endproof{$\hfill\Box$\\}
\def\R{\mathbb{R}}

\title[$W^{1,p}$-metrics and conformal metrics]{\bf $W^{1,p}$-metrics and conformal metrics with $L^{n/2}$-bounded
scalar curvature}

\author{Conghan Dong, Yuxiang Li, Ke Xu}
\address{Conghan Dong: Department of Mathematical Sciences\\Tsinghua University\\Beijing 100084\\P. R. China}
\email{dongch13@tsinghua.org.cn}
\address{Yuxiang Li: Department of Mathematical Sciences\\Tsinghua University\\Beijing 100084\\P. R. China}
\email{liyuxiang@tsinghua.edu.cn}
\address{Ke Xu: Department of Mathematical Sciences\\Tsinghua University\\Beijing 100084\\P. R. China}
\email{xuke16@mails.tsinghua.edu.cn}

\date{}
\begin{document}
\maketitle

\begin{abstract}
A $W^{1,p}$-metric on an $n$-dimensional closed Riemannian manifold naturally induces a distance function, provided $p$ is sufficiently close to $n$. If a sequence of metrics $g_k$ converges in $W^{1,p}$ to a limit metric $g$, then the corresponding distance functions $d_{g_k}$ subconverge to a limit distance function $d$, which satisfies $d\le d_g$.

As an application, we show that the above convergence result applies to a sequence of conformal metrics with $L^{n/2}$-bounded scalar curvatures, under certain geometric assumptions. In particular, in this special setting, the limit distance function $d$ actually coincides with $d_{g}$.

%In this paper we firstly show that  a $W^{1,p}$-metric
%on an $n$-dimensional smooth closed Riemannian manifold $(M,g_0)$ induces a distance function when $p$ is close to $n$.  If
%$g$ and $g^{-1}$, as tensors, are both in $W^{1,p}(M,g_0)$, the average limit of $g$ exists everywhere except an $(n-p)$-dimensional subset, so
%we can define a distance function $d_g$ by
%$$
%d_g(x,y)=\inf\left\{ \int_\gamma \sqrt{g(\dot{\gamma},\dot{\gamma})}: \mbox{piecewise smooth $\gamma$ from $x$ to $y$} \right\}.
%$$
%Let $\{g_k\}$ and $\{g_k^{-1}\}$ converge to $g$ and $g^{-1}$ respectively in $W^{1,p}(M,g_0)$. Then we prove that, after passing to a subsequence, $\{d_{g_k}\}$ converges to a distance function $d$, which in general satisfies  $d\leq d_g$.  As an application, we consider a sequence of conformal metrics $\{g_k=u_k^{\frac{4}{n-2}}g_0\}$ on $(M,g_0)$ with bounded $L^\frac{n}{2}$-norm of
%the scalar curvature. Under some geometric assumptions, we will show that, after passing to a subsequence, $\{g_k\}$ converges weakly to some $g=u^{\frac{4}{n-2}}g_0$ in $W^{2,q}(M,g_0)$ with $q<\frac{n}{2}$, and in this special case actually $\{d_{g_k}\}$  converges to $d_{g}$.
\end{abstract}

\section{Introduction}
%The problem considered in this paper can be stated as follows.
In this paper, we are interested in the convergence of a sequence of $W^{1,p}$-metrics on a Riemannian manifold, which is motivated by the study of conformal metrics with $L^{\frac{n}{2}}$-bounded scalar curvatures.

Let $(M, g_0)$ be a smooth closed $n$-dimensional Riemannian manifold and $p<n$ be sufficiently close to $n$. We first observe that, given a $W^{1,p}$-metric $g$ with respect to the background metric $g_0$, there is a well-defined distance function $d_g$ associated to $g$. Actually, using an idea similar to the Trace Embedding Theorem, we can show that $g$ is well-defined almost everywhere except a possible singular set of Hausdorff dimension at most $n-p$, which enables us to define
 $$
d_g(x,y):=\inf\left\{ \int_\gamma \sqrt{g(\dot{\gamma},\dot{\gamma})}: \mbox{piecewise smooth $\gamma$ from $x$ to $y$} \right\}, \forall x,y\in M.$$

Now suppose $\{g_k\}$ is a sequence of smooth metrics on $M$ such that $g_k$ and $g_k^{-1}$ converges to $g$ and $g^{-1}$ in $W^{1,p}(M,g_0)$, respectively. Then the limit $W^{1,p}$-metric $g$ induces a distance function $d_g$. On the other hand, the distance function $d_{g_k}$ associated to $g_k$ converges uniformly to a distance function $d$, as the metric space $(M,d_{g_k})$ converge to $(M,d)$ in sense of the Gromov-Hausdorff distance.
Then it is very natural to ask what is the relation between $d_g$ and $d$?

%To make that question precisely, we need to define an associated distance function from the limit metric $g$, where $g$ only lies in $W^{1,p}(M,g_0)$. The key idea is, given a function  $u\in W^{1,p}(U)$ with $U\subset \mathbb{R}^n$, we have canonical value of $u(x)$ away from a $(n-p)$-dimensional
%subset. Precisely, the set
%$$
%E:=\left\{x: \lim_{r\rightarrow 0}\frac{1}{|B_r(x)\cap U|}   \int_{B_r(x)\cap U}u\ \mbox{does not exist}\right\},
%$$
%has Hausdorff dimension at most $n-p$, and we say $u(x)$ is well-defined with value being
%its average limit centered at $x$ whenever $x\notin E$.
%Similarly in the Riemannian case, if $p$ is very close to $n$, $\sqrt{g(\dot\gamma,\dot\gamma)}$ is well-defined on a piecewise smooth curve.
%Then we can define
%$$
%d_g(x,y):=\inf\left\{ \int_\gamma \sqrt{g(\dot{\gamma},\dot{\gamma})}: \mbox{piecewise smooth $\gamma$ from $x$ to $y$} \right\}.$$
%It is not difficult to check that  $d\leq d_{g}$. However, we could only get $d=d_{g}$ when
%$g$ is continuous, which we believe is also true in general. In summary, our first main result
%is the following.

Our first main result is

\begin{theorem}\label{main1}
Let $(M, g_0)$ be an $n$-dimensional closed Riemannian manifold and $p\in(\frac{2n(n-1)}{2n-1},n]$. Suppose $\{g_k\}$ is a sequence  of smooth Riemannian metrics such that $\{g_k\} $ and $\{g^{-1}_k\}$ converge to $g$ and $g^{-1}$ in $W^{1,p}(M,g_0)$, respectively. Then, up to a subsequence, $\{d_{g_k}\}$ converges uniformly to a distance function $d$ with $d\leq d_g$. In particular, if $g$ and $g^{-1}$ are continuous, then $d=d_g$.
\end{theorem}

\begin{rem}
  Although here we only get the identity $d=d_g$ when $g$ is continuous, we believe it is still true without the assumption that $g$ and $g^{-1}$ are continuous.
\end{rem}

As an application, we next study the compactness of  a sequence of  conformal metrics $\{g_k=u_k^{\frac{4}{n-2}}g_0\}$ with $L^\frac{n}{2}$-bounded scalar curvature $\|R(g_k)\|_{L^{n/2}}\leq C$ and $\mathrm{Vol}(M,g_k)=1$.
Since $W^{2,\frac{n}{2}}$-space fails to be embedded into $C^0$, in general one can not expect point-wise compactness without additional assumptions, as shown by the counterexamples in \cite{Brendle,Brendle-Marques,Chang-Gursky-Wolff}. In fact, there is no compactness even if $R(g_k)$ is bounded $L^\infty$, see \cite{Brendle, Brendle-Marques}.

Here we assume in addition to the $L^\frac{n}{2}$-norm of the scalar curvature being bounded, its limiting measure is locally small as specified in the statement of the following theorem. We focus on the convergence of measures and  distance functions.  

\begin{theorem}\label{main2}
Let $(M,g_0)$ be an $n$-dimensional closed Riemannian manifold with $n\geq 3$. Suppose $\{g_k=u_k^{\frac{4}{n-2}}g_0\}$ is a sequence of conformal metrics such that $\mathrm{Vol}(M,g_k)=1$ and $|R(g_k)|^{\frac{n}{2}}dV_{g_k}$ converges weakly to a measure $\mu$ with $\mu(M)<\Lambda$.
Then, for any $q\in (1,\frac{n}{2})$, there exists $\varepsilon_0>0$, which only depends on $(M,g_0)$,  $\Lambda$ and $q$, such that
if
$$
\mu(\{x\})<\varepsilon_0,\forall x\in M,
$$
then, after passing to a subsequence, we have

1) $\{u_k\}$, $\{\frac{1}{u_k}\}$ and $\{\log u_k\}$ converge to $u$ , $\frac{1}{u}$ and $\log u$ weakly in $W_{loc}^{2,q}(M)$ respectively.

2)$d_{g_k}$ converges to $d_{g}$ in $C_{loc}^0(M\times M)$, where $d_g$ is the distance function associated to the limit metric $g:=u^\frac{4}{n-2}g_0.$
\end{theorem}

\begin{rem}
 Note that the limit metric $g\in W^{2,q}(M), q<\frac{n}{2}$.  $W^{2,q}$-space fails to be embedded into $C^0$. Hence in Theorem \ref{main2}, we can not get $g$ is continuous, but we can use the conformal condition to get the same conclusion as in Theorem \ref{main1}.
\end{rem}
\begin{rem}
Note that the $L^{\frac{n}{2}}$-norm of scalar curvature is rescaling invariant, i.e.
$$
\int_M|R(\lambda g_k)|^\frac{n}{2}dV_{\lambda g_k}=\int_M|R(g_k)|^\frac{n}{2}dV_{g_k}.
$$
So for a general sequence of collapsing metrics $\{g_k\}$, we can normalize the metric and set $g_k'=(c_ku_k)^\frac{4}{n-2}g_0$ such that $\mathrm{Vol}(M,g_k')=1$. Then after passing to a subsequence, the distance functions $d_{g_k'}$ associated to $g_k'$ converges uniformly to a distance function $d_g$, which is defined by $g=v^{\frac{4}{n-2}}g_0$, with $v$ being the $W^{2,q}$-weak limit of $c_ku_k$.
\end{rem}

In a recent paper \cite{Aldana-Carron-Tapie}, C. Aldana, G. Carron and S. Tapie obtained the above Gromov-Hausdoff convergence result in a similar setting, see \cite[Theorem 5.1]{Aldana-Carron-Tapie}). But they used a different method from ours, and the equality of $d$ and $d_g$ was left open.

In \cite{Li-Zhou}, the second author studied the bubble tree convergence of $\{g_k\}$ under a stronger assumption that  $\|R(g_k)\|_{L^p}<C$ with $p>\frac{n}{2}$. The compactness of conformal metrics with uniformly $L^p$-bounded sectional curvature was discussed in \cite{Chang-Yang1,Chang-Yang2,Gursky}.

%We should mention C. Aldana, G. Carron and S. Tapie \cite{Aldana-Carron-Tapie} showed that under some bounds on $L^{\frac{n}{2}}$-norm of scalar curvature, metric spaces induced by conformal Riemannian metrics are precompact in the Gromov-Hausdoff topology, and the limit is also a metric space defined by a conformal metric on $M$. Moreover, they also showed that there are $f\in W^{2,\frac{n}{2}}$ and $w\in C^\alpha$, $\alpha\in (0,1)$, such that $\log u=f+w$(see Theorem 5.1 in \cite{Aldana-Carron-Tapie}). But the method they used
%is different from ours. In our paper, Theorem \ref{main2} answers the question about the equality between $d$ and $d_g$, which is left open at the end of their paper. We also need to mention that, when $\int_M|K(g_k)|^{p}dV_{g_k}$ is bounded, where $p>\frac{n}{2}$ and $K(g_k)$ is the sectional curvature, the compactness has been discussed in
%\cite{Chang-Yang1,Chang-Yang2,Gursky}.

\vspace{1cm}

The rest of the paper is organized as follows. In Section 2, we study the basic properties of the average limit of a  $W^{1,p}$ function with $p<n$.
%Actually, those are some useful results about the Trace Embedding Theorem of $W^{1,p}$-functions in  \cite{Adams-Fournier,Evans-Gariepy,V.G.-T.O.}.
In Section 3, we discuss the convergence of a sequence of distance functions associated to $W^{1,p}$-metrics and prove Theorem~\ref{main1}. At last, we provide a key $\varepsilon$-regularity theorem and finish the proof of Theorem~\ref{main2} in Section~4.

%In section 4, we study metrics in a fixed conformal class. We use Moser iteration to prove $\varepsilon$-regularity at first. Then we prove our second main result. If there is no volume collapsing, then, up to a subsequence, $\{u_k\}$, $\{\log u_k\}$ and $\{u_k^{-1}\}$ converges to $u$, $\log u$ and $u^{-1}$ weakly in $W^{2,q}$ for some $q<\frac{n}{2}$. Moreover, a subsequence of distance functions $\{d_k\}$ uniformly converges to $d$. This gives the first conclusion in Theorem \ref{main2}. In the remaining parts, we use rescaling methods to show $d=d_g=u^\frac{4}{n-2}g_0$, even though $u$ is not continuous.

\textbf{Acknowledgment.}
 Part of this work was done while the third author was visiting S.-Y. A. Chang at Princeton University. She would like to thank S.-Y. A. Chang for  helpful discussions. The authors would like to thank the referees for their valuable suggestions on the revision.

\section{Traces of $W^{1,p}$-functions}
Recall that the \emph{Trace Embedding Theorem} states that the restriction of a $W^{1,p}$-function, which is defined on a domain in $\mathbb{R}^n$,  on a $k$-dimensional subset is a $L^q$-function, where $p<n$, $n-p<k\leq n$ and $p\leq q\leq\frac{np}{n-p}$,  see \cite[Theorem 4.12]{Adams-Fournier}. There is also the so-called \emph{trace inequality} in Sobolev Spaces, see \cite[Chapter 1]{V.G.-T.O.}.
Using similar ideas from those references, here we will show that the  \emph{centered average limit} of a given  $W^{1,p}$-function is well-defined almost everywhere, except on a subset of Hausdorff dimension smaller than $n-p$
(c.f. \cite{Evans-Gariepy}).

For $r>0$, we denote the $n$-ball centered at $x\in \mathbb{R}^n$ with radius $r$ by $B_r(x)$ and $B_r:=B_r(0)$. Given a function $u$ defined on a domain $U\subset \mathbb{R}^n$, the \emph{$r$-average} of $u$ at $x\in U$ is
$$
u_{x,r}:=\frac{1}{|B_r(x)\cap U|}\int_{B_r(x)\cap U}u(y)dy.
$$
%The \emph{average limit} of $u$ centered at $x$ is the limit $\lim_{r\to 0}u_{x,r}$, whenever it exists.

Now for $u\in W^{1,p}(B_2)$, define the singular set
$$
A(u)=\{x\in B_1:\lim_{\tau\rightarrow 0}osc_{r\in(0,\tau]} u_{x,r}>0\}.
$$
By Federer and Ziemer's theorem, we know $\dim A(u)\leq n-p$ (see \cite[Theorem 2.1.2]{Lin-Yang} or \cite[p.160]{Evans-Gariepy}).  Hence, we can define the \emph{centered average limit}
\begin{defi}
Given a function $u\in W^{1,p}(B_2), p<n$. For any $x\in B_1\setminus A(u)$, there exists $\hat{u}(x)$, such that
$$
\lim_{r\to 0}\frac{1}{B_r(x)}\int_{B_r(x)}|u(y)-\hat{u}(x)|dy=0,
$$
$\hat{u}(x)$ is called the \emph{centered average limit} of $u$ at $x$.
\end{defi}
Thus $\hat{u}$ is well-defined for $\mathcal{H}^s$-a.e. $x\in B_1$, with $s\in (n-p,n)$.

The following estimate will play an essential role in the next section.
\begin{lem}\label{measure.estimate}
Let $u\in W^{1,p}(B_2)$ and
$$
\mathcal{M}(u,t):=\{x\in B_1\setminus A(u): |\hat{u}|(x)>t\}.
$$
Assume $\|u\|_{L^1(B_2)}\leq \frac{t\omega_n}{4}$
and $s\in(n-p,n)$.
Then
$$
\mathcal{H}_{\infty}^s(\mathcal{M}(u,t))\leq\frac{\Lambda}{t^p}\int_{B_2}
|\nabla u|^p,
$$
where $\Lambda=\Lambda(n,s,p)$. Moreover, there exists
a cover $\{\overline{B_{r_i}(x_i)}\}$ of $\mathcal{M}(u,t)$, such that
$$
x_i\in\mathcal{M}(u,t),\s and\s \omega_s\sum_i r_i^s\leq\frac{\Lambda}{t^p}\int_{B_2}
|\nabla u|^p.
$$
\end{lem}
\proof Fix $x\in\mathcal{M}(u,t)$. By the definition of $\hat{u}$,
 $\frac{1}{|B_{r_0}|}|\int_{B_{r_0}(x)}u|>t$ for sufficiently small ${r_0}$.
We claim that there exists a small $r_0$ such that
\begin{equation}\label{lemma1-claim}
\frac{1}{{r_0}^s}\int_{B_{r_0}(x)}|\nabla u|^p\geq t_1:=(\frac{t}{\Lambda'})^p,
\end{equation}
where $\Lambda'=\Lambda'(n,p,s)$ is a constant to be determined later.

Assume for contradiction that the claim is not true. By the proof of Poincar\'e inequality (see \cite[pp.275-276]{Evans}), we have
$$
\frac{1}{|B_r|}\int_{B_r(x)}|u-u_{x,\frac{r}{2}}|^p\leq
\Lambda_1r^{p-n}\int_{B_r(x)}|\nabla u|^p,
$$
where $\Lambda_1$ only depends on $n$.
It follows
\begin{eqnarray*}
	|u_{x,r}-u_{x,\frac{r}{2}}|&=&\frac{1}{|B_r|}\left|\int_{B_r(x)}(u-u_{x,\frac{r}{2}})\right|\nonumber\\\nonumber
	&\leq&\frac{1}{|B_r|}\left(\int_{B_r(x)}|u-u_{x,\frac{r}{2}}|^p\right)^\frac{1}{p}|B_r|^{1-\frac{1}{p}}\\
	&=&\left(\frac{1}{|B_r|}\int_{B_r(x)}|u-u_{x,\frac{r}{2}}|^p\right)^\frac{1}{p}\nonumber\\
	&\leq&\left(\Lambda_1r^{p-n}\int_{B_r(x)}|\nabla u|^p\right)^\frac{1}{p}\nonumber\\
	&\leq&\Lambda_2r^\theta t_1^\frac{1}{p},
\end{eqnarray*}
where $\theta=\frac{p-n+s}{p}$ and $\Lambda_2=\Lambda_1^\frac{1}{p}$.

For $r_0\in[2^{-k},2^{-k+1})$, we have
$$
|u_{x,1}-u_{x,2^{-k}}|\leq \Lambda_2(\sum_{i=0}^{k-1} (2^{-i})^\theta)t_1^\frac{1}{p}
\leq \Lambda_3 t_1^\frac{1}{p},
$$
and
\begin{eqnarray*}
	|u_{x,2^{-k}}-u_{x,r_0}|&=&\frac{1}{|B_{r_0}|}\left|\int_{B_{r_0}(x)}(u-u_{x,2^{-k}})\right|\\
	&\leq&\frac{|B_{2^{-k+1}}|}{|B_{r_0}|}\frac{1}{|B_{2^{-k+1}}|}\int_{B_{2^{-k+1}(x)}}\left|u-u_{x,2^{-k}}\right|\\
	&\leq&2^n\Lambda_2(2^{1-k})^\theta t_1^\frac{1}{p},
\end{eqnarray*}
where $\Lambda_3=2\Lambda_2$.
Then
\begin{eqnarray}\label{decay}
|u_{x,1}-u_{x,r_0}|\leq\Lambda_3t_1^\frac{1}{p}=
\frac{\Lambda_3}{\Lambda'}t.
\end{eqnarray}
Note that $|u_{x,1}|\leq\frac{t}{4}$, so we get a
contradiction if we set $\Lambda'>2\Lambda_3$. This proves our claim~(\ref{lemma1-claim}).

To complete the proof of the lemma, note that by  Vitali Covering Theorem,  there exist
pairwise disjoint closed balls
$\{\overline{B_{r_i}(x_i)}\}_{i=1}^\infty$ such that
$$
\frac{1}{r_i^s}\int_{B_{r_i}(x_i)}|\nabla u|^p\geq t_1,\s
\mathcal{M}(u,t)\subset\bigcup_i\overline{B_{5r_i}(x_i)}.
$$
Therefore, we get
\begin{eqnarray*}
\mathcal{H}_{\infty}^s(\mathcal{M}(u,t))&\leq& \sum_i\omega_s(5r_i)^s=5^s\omega_s\sum_i r_i^s\\
&\leq&
\frac{1}{t_1}5^s\omega_s\int_{\cup B_{r_i}(x_i)}|\nabla u|^p\\
&\leq&\frac{1}{t_1}5^s\omega_s\int_{B_2}|\nabla u|^p.
\end{eqnarray*}
\endproof

As an application of Lemma~\ref{measure.estimate}, we show that $W^{1,p}$-convergence implies
$\mathcal{H}^s$-a.e. convergence for each $s>n-p$.
\begin{lem}\label{equal}
Assume $u_k, u\in W^{1,p}(B_2)$ and $\|u_k-u\|^p_{W^{1,p}(B_2)}<\frac{1}{2^k}$, then for any   $s>n-p$, $\hat{u}_k$ converges to
$\hat{u}$ for $\mathcal{H}^{s}$-a.e. $x\in B_1$.
\end{lem}
\proof
Set
$$
A=(\bigcup_{i=1}^\infty A(u_k))\bigcup A(u),\s
E_{km}=\{x\in B_1\setminus A:|\hat{u}_k-\hat{u}|<\frac{1}{m}\},
$$
and
$$
E=\bigcap_{m=1}^\infty\bigcup_{i=1}^\infty\bigcap_{k=i}^\infty
E_{km}.
$$
It is easy to check that for any $x\in E$,  $\hat{u}_k(x)$ converges to $\hat{u}(x)$.

Let
$$
F=E^c\cap B_1\setminus A=\bigcup_{m=1}^\infty\bigcap_{i=1}^\infty\bigcup_{k=i}^\infty
F_{km},
$$
where
$$
F_{km}=\{x\in B_1\setminus A:|\hat{u}_k-\hat{u}|\geq\frac{1}{m}\}.
$$
Since $\hat{u}_k-\hat u=\widehat{u_k-u}$ and $\frac{1}{|B_1(x)|}
\int_{B_1(x)}|u_k-u|dx\rightarrow 0$,
by Lemma \ref{measure.estimate},
$\mathcal{H}^s_{\infty}(F_{km})\leq Cm2^{-k}$ when $k$
is sufficiently large. It follows that
$$
\mathcal{H}_{\infty}^s(\bigcap_{i=1}^\infty\bigcup_{k=i}^\infty
F_{km})=0,
$$
which implies $\mathcal{H}^s(F)=0$. Since $B_1\setminus E
\subset A\cup F$, we get $\mathcal{H}^s(B_1\setminus E)=0$.
\endproof

\begin{rem}\label{second}
Lemma \ref{equal}  provides another approach to define the value of $u$ at a point. Select a sequence of smooth functions
$u_k$ satisfying  $\|u_k-u\|_{W^{1,p}}<2^{-k}$.
Since $\hat{u}_k=u_k$, by Lemma \ref{equal},
$u_k$ converges to $\hat{u}$ for $\mathcal{H}^{s}$-a.e. $x$ whenever $s>n-p$. Therefore, $\hat{u}$ is in fact an $\mathcal{H}^s$-a.e. limit of $u_k$. Using this point of view, one can easily check the following:

 1) when $f\in C^1$, $\widehat{f u}=f\widehat{u}$ for $\mathcal{H}^s$-a.e. $x$.

 2) when $s>n-\frac{p}{q}$, $\widehat{u^q}=\hat{u}^q$ for $\mathcal{H}^s$-a.e. $x$.
\end{rem}

Let $p>n-m$ and $\Sigma$ be a  compact $m$-dimensional submanifold of $B_1$.  We can establish a trace embedding inequality  on $\Sigma$. Applying Theorem 1.1.2
in \cite{V.G.-T.O.}
to $\mu=\mathcal{H}^m\lfloor\Sigma$, we have
$$
\|u\|_{L^1(\Sigma)}\leq C(\Sigma)\|u\|_{W^{1,m}(\R^n)},
$$
where $u\in C_0^\infty(\R^n)$.
Given a function $u\in W^{1,p}(B_1)$, after extending it to a function $u'\in W^{1,p}_{0}(B_2)$ with
$$
\|u'\|_{W^{1,p}(B_2)}\leq C(n) \|u\|_{W^{1,p}(B_1)},
$$
we can find $u_k\in C^\infty_0(B_2)$ such that $\|u_k-u'\|_{W^{1,p}(\R^n)}\rightarrow 0$.
Then $\{u_k\}$ is a Cauchy sequence in
$L^1(\Sigma)$.
By Lemma \ref{equal}, we may assume
$u_k$ converges to $\hat{u}$ for $\mathcal{H}^m-$a.e. $x\in\Sigma$. Therefore, we obtain
\begin{equation}\label{trace.inequality}
\int_\Sigma|\hat{u}|d\mathcal{H}^m\lfloor\Sigma\leq C(\Sigma)
\|u\|_{W^{1,p}(B_1)}. \\
\end{equation}

Now we consider the case when $(M,g)$ is a smooth Riemannian manifold and $u\in W^{1,p}(M)$. For $x\in M$, in a local coordinate chart $(x^1,\cdots,x^n)$ centered at $x$, we can define $\hat{u}(x)$ to be the limit of $\frac{1}{\omega_n r^n}
\int_{B_r}udx$ as $r\rightarrow 0$, where $B_r$ is the Euclidean ball as before. In view of Remark \ref{second}, one checks that the value of $\hat{u}(x)$ is independent of the choice of coordinate chart for $\mathcal{H}^s-$ a.e.
$x\in M$, where $s>n-p$. 
%Moreover, we have the following:
\begin{lem}\label{mf}
There exists a subset $E\subset M$ with dimension smaller than $(n-p)$ such that
for any $x\notin E$, there holds
\begin{equation}\label{mfd}
\hat{u}(x)=\lim_{r\to 0}\frac{1}{\mathrm{Vol}(B^{g}_r(x))}\int_{B^{g}_r(x)} udV_g,
\end{equation}
where
$$
B_r^g(p)=\{x\in M : d_g(x,p)<r\}.
$$
\end{lem}
\proof
Locally in a coordinate $(x^1,\cdots,x^n)$, we set
$$
\Lambda_s=\{x:\varlimsup_{r\to 0} \frac{1}{r^s}\int_{B_r^n(x)}|\nabla u|^pdx>0\},\s
and \s
\Lambda=\bigcap_{s\in(n-p,n)}\Lambda_s.
$$
It is a standard result that $H^s(\Lambda_s)=0$ when $s\in(n-p,n)$ (c.f. \cite[Lemma 2.1.1]{Lin-Yang}). Since $\Lambda_{s'}\subset\Lambda_{s}$ for any $s'<s$, we have $\dim \Lambda<n-p$. We will
show that \eqref{mfd} holds for any $x\notin \Lambda$. Obviously,
we only need to prove \eqref{mfd} holds
for any $x\notin \Lambda_s$ and $s\in (n-p,n)$.

Fix an $x_0\notin \Lambda$. As in
\eqref{decay}, we have
$$
|u_{x,r}-u_{x,r'}|\leq \Lambda_3\left(\frac{1}{r^s}\int_{B_r(x)}|\nabla u|^pdx\right)^\frac{1}{p},\s whenever\s r'<r.
$$
Thus $u_{x,r}$ converges as $r\rightarrow 0$  for any $x\notin \Lambda_s$. Denoting
 $u_r(x)=u(x_0+rx)$ and applying the Poincar\'e inequality for a ball with any fixed radius $R>0$, we get
$$
\int_{B_R}\left|u_r(x)-\frac{1}{|B_1|}\int_{B_1}u_rdx\right|dx\leq
\int_{B_R}|\nabla u_r|^pdx=R^sr^{s+p-n}\frac{1}{(Rr)^s}\int_{B_{Rr}(x_0)}|\nabla u|^pdx
\rightarrow 0.
$$
Since $\frac{1}{|B_1|}\int_{B_1}u_rdx=u_{x_0,r}$ converges to $\hat{u}(x_0)$, we have
$$
\lim_{r\rightarrow 0}\int_{B_R}\left|u_r(x)-\hat{u}(x_0)\right|dx=0.
$$
Note that, when $r$ is sufficiently small, we may assume $B_1^{g(x_0+rx)/r^2}(0)\subset B_R$. Then we conclude
\begin{eqnarray*}
& \ &\lim_{r\to0}\frac{1}{\mathrm{Vol}(B_r^{{g}}(x_0))}\int_{B_{r}(x_0)} |u-\hat{u}(x_0)|dV_g\\
&=&\lim_{r\to0}\frac{1}{\mathrm{Vol}(B_1^{{g(x_0+rx)}/r^2}(0))}\int_{B_1^{{g(x_0+rx)}/r^2}(0))} |u_r-\hat{u}(x_0)|dV_{g(x_0+rx)/r^2}\\
&\leq& C(M)\cdot\lim_{r\to0}\int_{B_R}|u_r-\hat{u}(x_0)|dx\\
&=&0.
\end{eqnarray*}
\endproof

\section{$W^{1,p}$-metrics}
Suppose $(M,g_0)$ is a smooth closed $n$-manifold.
Let $g$ be a symmetric  tensor of type $(0,2)$, which is positive almost everywhere.
Let $g^{-1}$ be the corresponding inverse tensor of type $(2,0)$.
We say that $g$ is a $W^{1,p}$-metric if both $g$ and $g^{-1}\in
W^{1,p}_{loc}(M,g_0)$. The goal of this section is to define the distance functions induced by $W^{1,p}$-metrics and study the compactness of such metrics.  In a local coordinate chart, we can write
$$
g=g_{ij}dx^i\otimes dx^j,\s  and\s
g^{-1}=g^{ij}\frac{\partial}{\partial x^i}\otimes \frac{\partial}{\partial x^j}.
$$
Then the functions $g_{ij}$, $g^{ij}$ belong to $W^{1,p}_{loc}$, and $(g^{ij})(g_{ij})=I$ as matrices.
Now define the centered average limit of the metric by
$$
\hat{g}_{ij}(x)= \lim_{r\to 0}\frac{1}{|B_r(x)|}\int_{B_r(x)} g_{ij}(y)dy,
$$
and the corresponding tensor by
$$
\hat{g}(x)(V(x),V(x))=\sum_{i,j}\hat{g}_{ij}(x)V_i(x)V_j(x), \quad \forall V\in \Gamma(TM).
$$
By Remark \ref{second}, when $s>n-p$, $\hat{g}$ and $\widehat{g^{-1}}$ are well-defined on
$T_xM$ and $(T_xM)^*$ for $\mathcal{H}^s$-a.e. $x$. Moreover, by Lemma \ref{mf}
\begin{align*}
\hat{g}(x)(V(x),V(x))&=\lim_{r\to 0}\frac{1}{\mathrm{Vol}(B_r^{g_0}(x))}\int_{B_r^{g_0}(x)}g(y)(V(y),V(y))dV_{g_0}
\\ &= \sum_{i,j} \left(\lim_{r\to 0}\frac{1}{|B_r(x)|}\int_{B_r(x)} g_{ij}(y)dy\right)V_i(x)V_j(x)
\end{align*}

Next, we define the associated distance function by
$$
d_g(x,y)=\inf\left\{ \int_\gamma \sqrt{g(\dot{\gamma},\dot{\gamma})}: \mbox{piecewise smooth $\gamma$ from $x$ to $y$} \right\}.$$
First of all, we need to show that $d_g$ is indeed a distance function.
\begin{lem}
When $p\in (n-1,n)$,
$d_g$ is a distance function and it's continuous on $M\times M$.
\end{lem}

\proof
We first show that $d_g(x,y)<+\infty$ for any $x, y\in M$.

Let $\varphi_x$ be the exponential map from $T_xM$
to $M$.
Since $(M,g_0)$ is compact, there exists a number $\tau=\tau(M,g_0)>0$, such that for each $x\in M$, $\varphi_x$ induces  normal coordinates $({x'}^1,\cdots,{x'}^n)$ on $B_\tau^{g_0}(x)$
with
$$
|g_{0,ij}(x')-\delta_{ij}|<\frac{1}{2}.
$$
It follows that the metric
$$
g=g_{ij}(x')d{x'}^i\otimes d{x'}^j,
$$
satisfies
\begin{equation}\label{wip.chart}
\frac{1}{C}\|g\|_{W^{1,p}(B_\tau^{g_0}(x),g_0)}\leq\|(g_{ij})\|_{W^{1,p}(B_\tau^n(0))}\leq C\|g\|_{W^{1,p}(B_\tau^n(0))},
\end{equation}
and
\begin{equation}\label{wip.chart2}
\frac{1}{C}\|g^{-1}\|_{W^{1,p}(B_\tau^{g_0}(x),g_0)}\leq\|(g^{ij})\|_{W^{1,p}(B_\tau^n(0))}\leq C\|(g^{-1})\|_{W^{1,p}(B_\tau^n(0))},
\end{equation}
where $C$ is a constant independent of $x$. Following \cite[p.178]{Petersen}, we call a curve $\gamma_{pq}$ joining $p$ and $q$ a \emph{segment}in $(M,g_0)$ if $Length(\gamma_{pq})=d_{g_0}(p,q)$ and $|\dot\gamma_{pq}|$ is constant.
Let $\gamma:[0,l]\rightarrow M$ be the segment in $(M,g_0)$
from $x$ to $y$. Then we can find $x_1=\gamma(t_1)$, $x_2=\gamma(t_2)$, $\cdots$,
$x_m=\gamma(t_m)$, such that
$$
m<\frac{2l}{\tau}, \s 0\leq t_{i+1}-t_i\leq \frac{\tau}{2}.
$$
For convenience, we set $x_0=\gamma(0)=x$ and $x_{m+1}=\gamma(t_{m+1})=y$.

It suffices to prove $d_{g}(x_i,x_{i+1})<+\infty$.
Without loss of generality, we assume
the coordinate of $x_{i+1}$ in  chart
$(B_\tau^{g_0}(x_i),\varphi_{x_i}^{-1})$ is
$(\delta_i,0,\cdots,0)$, where
$\delta_i=t_{i+1}-t_i$.
Obviously,
$$
d_g(x_i,x_{i+1})\leq\int_0^{\delta_i}\sqrt{\hat{g}_{11}(t,0,\cdots,0)}dt\leq
\delta_i^\frac{1}{2}\sqrt{\int_0^{\delta_i} \hat{g}_{11}(t,0,\cdots,0)dt}.
$$
By \eqref{wip.chart}, \eqref{wip.chart2} and \eqref{trace.inequality},
$\hat{g}_{11}(t,0,\cdots,0)$ is integrable on $[0,\delta_i]$. It follows that $d_g(x,y)<+\infty$ and $d_g(x,y)\rightarrow 0$ when $l\to 0$.

Next, we prove that $d_g(x,y)>0$ for any $x\neq y$. In fact, we can
prove a stronger result here: for any $\delta>0$, there exists $\delta'>0$, which depends on $g_0$, $\delta$ and
$||g^{-1}||_{W^{1,p}(M,g_0)}$, such that if $d_{g_0}(x,y)\geq\delta$,
then
\begin{equation}\label{positive.distance}
d_g(x,y)\geq \delta'.
\end{equation}

Assume $\gamma:[0,l]\rightarrow M$ is an arbitrary piecewise smooth curve from $x$ to $y$ in $M$. We set
$$
l'=\sup\{t\in[0,l]:\gamma([0,t])\subset B_\frac{\tau}{2}^{g_0}(x)\}.
$$
Obviously,
$$
d_{g_0}(x,\gamma(l'))=\min\{\tau/2,d_{g_0}(x,y)\}.
$$
In the coordinate defined by $\varphi_x$, we set
$\|g_{ij}\|=\sqrt{\sum (g_{ij})^2}$.
It is well-known that
$\|(g_{ij})\|^2$ is the
quadratic sum  of the eigenvalues of $(g_{ij})$. Let $\lambda$ be the
smallest eigenvalue of $(\hat{g}_{ij})$. Since $\frac{1}{\lambda}$ is an
eigenvalue of $(\hat{g}^{ij})$, we have
$$
E_a:=\{x'\in B_\frac{\tau}{2}^n(0):\lambda(x')<a\}\subset\{x'\in B_\frac{\tau}{2}^n(0):\|(\hat{g}^{ij})\|(x')>\frac{1}{a}\}.
$$
By the inequality $\|(\hat{g}^{ij})\|\leq c(n)\sum_{ij} |\hat{g}^{ij}|$, together with
Lemma \ref{measure.estimate},  we can find a sufficiently small $a$, which depends on $||g^{-1}||_{W^{1,p}(M)}$, $\delta$ and $\tau$, such that
$$
\mathcal{H}_{\infty}^1(\varphi^{-1}_x(\gamma|_{[0,l']})\cap E_a)\leq\mathcal{H}_\infty^1(\{x'\in B_\frac{\tau}{2}^n(0):\|\hat{g}^{ij}\|>\frac{1}{a}\})\leq \sum_{ij}a^{p}\Lambda'\|{g}^{ij}\|^p_{W^{1,p}(B_\frac{\tau}{2}^n(0))}< \frac{d_{g_0}(x,\gamma(l'))}{4}.
$$
Note that using balls to cover $\varphi^{-1}_x(\gamma|_{[0,l']})$ might increase $\mathcal{H}_{\infty}^1(\varphi^{-1}_x(\gamma|_{[0,l']}))$ by at most a factor of $2$ (see \cite{Han-Lin}). It follows that
$$
d_{\R^n}(0,\varphi_x^{-1}(\gamma(l')))=d_{g_0}(x,\gamma(l')) \leq 2\mathcal{H}_{\infty}^1(\varphi^{-1}_x(\gamma|_{[0,l']})),
$$
which implies
\begin{eqnarray*}
\mathcal{H}^1(\varphi^{-1}_x(\gamma|_{[0,l']})\setminus E_a)&\geq&\mathcal{H}_{\infty}^1(\varphi^{-1}_x(\gamma|_{[0,l']})\setminus E_a) \geq\mathcal{H}^1_{\infty}(\varphi^{-1}_x(\gamma|_{[0,l']}))-\mathcal{H}^1_{\infty}(\varphi^{-1}_x(\gamma|_{[0,l']})\cap E_a)\\
&\geq&\frac{d_{g_0}(x,\gamma(l'))}{2} -\frac{d_{g_0}(x,\gamma(l'))}{4}=\frac{1}{4}d_{g_0}(x,\gamma(l')).
\end{eqnarray*}
Since $\gamma$ is locally Lipschitz continuous,
$$
\int_{\varphi_x^{-1}(\gamma|_{[0,l']}) \setminus E_a} |\dot\gamma|\geq \mathcal{H}^1(\varphi^{-1}_x(\gamma|_{[0,l']})\setminus E_a).
$$
we have
$$
\int_\gamma\sqrt{\hat{g}(\gamma)(\dot\gamma,\dot\gamma)}
\geq\int_{\varphi_x^{-1}(\gamma|_{[0,l']})\setminus E_a}\sqrt{\hat{g}_{ij}(\gamma)\dot\gamma^i\dot\gamma^j}
\geq\int_{\varphi_x^{-1}(\gamma|_{[0,l']}) \setminus E_a}\lambda|\dot\gamma|\geq\frac{a}{4}d_{g_0}(x,\gamma(l')).
$$
This completes the proof, by letting $\delta'=\frac{a}{4}\min\{\tau/2,\delta\}$.

\endproof

{\it The proof of Theorem \ref{main1}:}
Since $|\nabla_{g_k,x} d_{g_k}(x,y)|=1$, in  local coordinates we have
$$
\lambda(x)|\nabla_{x} d_{g_k}(x,y)|\leq 1,
$$
where $\lambda$ is the smallest eigenvalue of $(g_{k,ij})$.
Since $\frac{1}{\lambda}$ is an eigenvalue of
$g_k^{ij}$, we see
$$
|\nabla_x d_{g_k}(x,y)|\leq\frac{1}{\lambda}\leq c(n)\sum_{ij} |g_k^{ij}(x)|.
$$
Similarly, we have $|\nabla_y d_{g_k}(x,y)|<c(n)\sum_{ij}|g_k^{ij}(y)|$,
hence $d_{g_k}$ is bounded in $W^{1,\frac{np}{n-p}}(M\times M,g_0)$. Therefore
we may assume $d_{g_k}$ converges to a function $d$ in $C^0(M\times M)$.
By \eqref{positive.distance}, we can assume further
\begin{equation}\label{positive.distance2}
d(x,y)\geq \tau,\s \mbox{whenever} \s d_{g_0}(x,y)>\delta.
\end{equation}

Next, we  prove that $d\leq d_{g}$. Given two points $x,y \in M$, take a piecewise smooth curve $\gamma$ from $x$ to $y$.
By \eqref{trace.inequality}, after passing to a subsequence,  $\sqrt{g_{k}(\dot{\gamma},\dot{\gamma}})$
converges to $\sqrt{g(\dot{\gamma},\dot{\gamma})}$
for $\mathcal{H}^1$-a.e. $x\in \gamma$. Moreover, we have
$$
\int_\gamma\left(\sqrt{g_k(\dot{\gamma},\dot{\gamma})}\right)^2<C,
$$
which implies
$$
\lim_{k\rightarrow+\infty}\int_\gamma \sqrt{g_k(\dot{\gamma},\dot{\gamma})}=\int_\gamma\sqrt{g(\dot{\gamma},\dot{\gamma})}.
$$
Thus, we arrive at $d(x,y)\leq d_g(x,y)$.

Finally, we show that $d=d_{g}$ in the case when
$g$ and $g^{-1}$ are continuous. For any $\varepsilon>0$ fixed, let
$$
E_k=\{x:g_k>(1-\varepsilon)g\}=\{x:g-g_k<\varepsilon g\},\s and\s F_k=E_k^c.
$$
We claim that
$$
\lim_{k\rightarrow+\infty}\mathcal{H}^1_\infty(F_k)=0.
$$
Since $M$ is compact, we only need to prove the claim in a local coordinate chart $\varphi: U\rightarrow \R^n$. That is, we only need to check that for any $B_R \subset \varphi(U)$,
$$
\mathcal{H}^1_\infty(B_R\cap F_k)\rightarrow 0.
$$

For simplicity,  we denote the
maximum eigenvalue and the minimum eigenvalue of a matrix $A$ by
$\Lambda(A)$ and $\lambda(A)$ respectively.
Since $g$ and $g^{-1}$ are continuous, we assume for any
$x\in B_R$,
$$
\frac{\lambda(g_{ij}(x))}{\|(g_{ij}(x))\|}\geq\varepsilon_1
$$
for some $\varepsilon_1>0$. Note that
\begin{eqnarray*}
F_k\cap B_R&\subset&\{x\in B_R:\Lambda(g_{ij}-g_{k,ij})\geq\varepsilon\lambda(g_{ij})\}\\
&\subset&\{x\in B_R:\|g_{ij}-g_{k,ij}\|\geq\varepsilon\lambda(g_{ij})\}\\
&\subset&\{x\in B_R:\|g_{ij}\|\cdot\|I-g_{k,ij}g^{ij}\|\geq\varepsilon\lambda(g_{ij})\}\\
&\subset&\{x\in B_R:\|I-g_{k,ij}g^{ij}\|\geq\varepsilon_1\epsilon\}.
\end{eqnarray*}

From the identity
$$
\nabla(g_{k,ij})(g^{ij})=\nabla (g_{k,ij})(g^{ij})+
(g_{k,ij})(\nabla g^{ij}),
$$
we see
$$
\|I-(g_{k,ij})(g^{ij})\|_{W^{1,q}}\rightarrow 0,
$$
for any $q<\frac{np}{2n-p}$. Since $p>2n\frac{n-1}{2n-1}$, we can
choose $q$ such that $q>n-1$.
By Lemma \ref{measure.estimate}, after passing to a subsequence, we get
$$
\lim_{k\rightarrow+\infty}\mathcal{H}^1_\infty(\{x\in B_R:\|I-(g_{k,ij})(g^{ij})\|\geq\varepsilon_1\varepsilon\})=0.
$$
Thus  the claim follows.

The above claim implies that, given $\varepsilon'>0$, for any $k$ sufficiently large, we can cover $F_k$, which is a compact subset,
by finitely many
balls $\overline{B_{r_1}}(x_1)$, $\cdots$, $\overline{B_{r_m}}(x_m)$, such that
$$
\sum r_i<\varepsilon'.
$$
Let $C_1$, $\cdots$, $C_{m'}$ be the connected components of
$B=\bigcup \overline{B_{r_i}(x_i)}$ and
set $t_1=\inf\{t:\gamma(t)\in B\}$. Without loss of generality, we assume $\gamma(t_1)\in C_1$. Put $t_2=\sup\{t:\gamma(t)\in C_1\}$, and replace
$\gamma|_{[t_1,t_2]}$ with the segment $\overline{\gamma(t_1)\gamma(t_2)}$.
In the same manner, we may choose $t_3=\inf\{t:\gamma(t)\in B\setminus C_1\}$  and by induction,
we can find
$$
0\leq t_1<t_2<t_3<\cdots<t_{m'}\leq 1,
$$
such that
$$
\sum_id_{g_0}(\gamma(t_{2i}),\gamma(t_{2i-1}))\leq\sum_idiam(C_i)\leq\sum r_i<\varepsilon'.
$$
This give rise to a new curve $\gamma'$ in place of $\gamma$, such that
Then
\begin{align*}
\int_{\gamma} \sqrt{g_k(\dot\gamma,\dot\gamma)} &\geq \int_{\gamma\cap\gamma'}
\sqrt{g_k(\dot\gamma,\dot\gamma)}\\
&\geq
(1-\varepsilon)^\frac{1}{2}(
\int_{\gamma'}\sqrt{g({\dot{\gamma}}',{\dot{\gamma}}')}-\sum_i\int_{\gamma(t_{2i-1})}^{\gamma(t_{2i})}\sqrt{g(\dot{\gamma'},\dot{\gamma'}}))
\\
&\geq
(1-\varepsilon)^\frac{1}{2}(d_g(x,y)-\varepsilon'\|\sqrt{g}\|_{C^0}).
\end{align*}
Therefore,
$$
d_{g_k}(x,y)\geq (1-\varepsilon)^\frac{1}{2}(d_g(x,y)-\varepsilon'\|\sqrt{g}\|_{C^0}).
$$
Now let $k\rightarrow+\infty$, then $\varepsilon'\rightarrow 0$,
and  finally $\varepsilon\rightarrow 0$, we get
$$
d(x,y)\geq d_g(x,y).
$$
\endproof

Although we only consider a compact manifold in Theorem \ref{main1}, a similar result
actually holds for certain complete manifolds. For example, we have
\begin{pro}\label{complete.case}
Let $\{g_k\}$ be a sequence of metrics defined on $\R^n$ and
assume  $g_k$ and $g_k^{-1}$ converge to $g_{\R^n}$ and $g_{\R^n}^{-1}$ respectively
in $W_{loc}^{1,p}(\R^n)$ for some $p>2n\frac{n-1}{2n-1}$. Then, after passing to a subsequence, $d_{g_k}(x,y)$ converges to
$|x-y|$.
\end{pro}
\proof
Let $R=|x-y|$. We only need to prove this proposition on $\overline{B_{2R}}$. We omit the details since the proof is almost the same with the one of Theorem \ref{main1}.
\endproof

\section{Conformal metrics with  $L^\frac{n}{2}$-bounded scalar curvature}\label{conformal.metrics}
First we recall some notations in conformal geometry. Let $(M,g)$ be a closed Riemannian manifold. Denote the  scalar curvature by $R(g)$ (or $R_g$). Let $g=u^{\frac{4}{n-2}}g_0$ be a conformal metric, then $u$ satisfies the following equation
$$
-\frac{4(n-1)}{n-2}\Delta u+R(g_0)u=R(g)u^\frac{n+2}{n-2}.
$$

\subsection{$\varepsilon$-regularity}
Again we denote by $B_r$ a ball in $\R^n$ with radius $r$, centered at $0$.  Let  $u$ be a weak solution of
\begin{equation}\label{equation.epsilon}
-div(a^{ij}u_{j})=fu,
\end{equation}
where
\begin{equation}\label{aij}
0<\lambda_1 I\leq (a^{ij}),\s \|a^{ij}\|_{C^0(B_2)}+\|\nabla a^{ij}\|_{C^0(B_2)}
<\lambda_2.
\end{equation}
%Then we have the following:
\begin{lem}\label{Lalpha}
Suppose $u\in W^{1,2}(B_2)$ is a positive weak solution of equations \eqref{equation.epsilon} and  \eqref{aij}. Assume
$$
\int_{B_2}|f|^\frac{n}{2}\leq \Lambda,
$$
then
$$
r^{2-n}\int_{B_r(x)}|\nabla\log u|^2<C,\s \forall B_r(x)\subset B_1.
$$
Moreover, there exist constants $\alpha$ and $C$, which depend on $\Lambda$, $\lambda_1$, $\lambda_2$,  such that
$$
\int_{B_1}(cu)^{\alpha}+\int_{B_1}(cu)^{-\alpha}<C,
$$
where $-\log c$ is the mean value of $\log u$ on $B_1$.
\end{lem}
\proof
 For a ball $B_{2r}(x)\subset B_2(0)$, take $\phi=\eta^2u^{-1}$ as a test function, with $\eta\equiv 1$ on $B_r(x)$, $\eta\in C^\infty_0(B_{2r}(x))$ and $|\nabla \eta|\leq \frac{C}{r}$. Multiplying
\eqref{equation.epsilon} by $\phi$ and integrating, we get
$$
\int_{B_{2r}(x)}\eta^2u^{-2}|\nabla u|^2\leq C\left(\int_{B_{2r}(x)} |\nabla\eta|^2+(\int_{B_{2r}(x)} f^{\frac{n}{2}})^{\frac{2}{n}}(\int_{B_{2r}(x)} \eta^{\frac{2n}{n-2}})^{\frac{n-2}{n}}\right),
$$
which implies
\begin{align*}
\int_{B_{r}(x)}|\nabla \log u|^2&\leq  Kr^{n-2}.
\end{align*}
By the Sobolev Embedding Theorem and the John-Nirenberg Lemma \cite[Theorem 3.5]{Han-Lin}, for $\alpha=\frac{C(n)}{K}$, we have
\begin{equation}\label{JN}
\|u\|_{L^{\alpha}(B_{1})}\|u^{-1}\|_{L^{\alpha}(B_1)}\leq C.
\end{equation}
Let $v=\log cu$, where $c$ is chosen such that
$$
\int_{B_1}v=0.
$$
By the Poincar\'e inequality, we can assume
$$
\int_{B_1}|v|\leq \beta_0,
$$
where $\beta_0$ only depends on $\Lambda$, $\lambda_1$ and
$\lambda_2$. Denote the Lebesgue measure over $\R^n$ by $L^n$.
Let
$$
E=\{x:v\leq \frac{2\beta_0}{L^n(B_1)}\}.
$$
Then
$$
L^n(B_1\setminus E)\leq \frac{L^n(B_1)}{2\beta_0}\int_{B_1}|v|\leq\frac{L^n(B_1)}{2},
$$
and $L^n(E)\geq\frac{1}{2}L^n(B_1)$.
By \eqref{JN}, we get
$$
C\geq \int_{B_1}(cu)^\alpha\int_{B_1}(cu)^{-\alpha}\geq\int_{B_1}(cu)^{\alpha}\int_{E}(cu)^{-\alpha}\geq \frac{1}{2}L^n(B_1)
e^{-\frac{2\alpha \beta_0}{L^n(B_1)}}\int_{B_1}(cu)^{\alpha}.
$$

In the same way, we  can get the estimate of $\int_{B_1}(cu)^{-\alpha}$.
\endproof

\begin{lem}\label{regularity}
Suppose $u\in W^{1,2}(B_2)$ is a positive solution of \eqref{equation.epsilon}, \eqref{aij}, and $\log u\in W^{1,2}(B_2)$.
Then for any $q\in (0,\frac{n}{2})$, there exists $\varepsilon_0
=\varepsilon(q,\lambda_1,\lambda_2)>0$, such that
if
$$
\int_{B_2}|f|^\frac{n}{2}<\varepsilon_0,
$$
then
$$
\|\nabla\log u\|_{W^{1,q}(B_{\frac{1}{2}})}\leq C(\lambda_1,\lambda_2,\epsilon_0).
$$
and
$$
e^{-\frac{1}{|B_\frac{1}{2}|}\int_{B_\frac{1}{2}}\log u}\|u\|_{W^{2,q}(B_\frac{1}{2})}+e^{\frac{1}{|B_\frac{1}{2}|}\int_{B_\frac{1}{2}}\log u}\|u^{-1}\|_{W^{2,q}(B_\frac{1}{2})}
\leq C(\lambda_1,\lambda_2,\epsilon_0).
$$
\end{lem}
\proof  Let $v=\log u$. In order to apply Lemma \ref{Lalpha}, we firstly assume $\int_{B_1}v=0$.

Let $\eta$ be a smooth cutoff function and $\phi=\eta^2u^{\beta}$ be a test function, where $\eta$ and $\beta\neq -1$ or $0$ will be defined later. Multiplying both sides of
(\ref{equation.epsilon}) by $\phi$ and integrating, we obtain
$$
\int_{B_1} 2\eta\nabla\eta u^{\beta}\nabla u+\int_{B_1} \eta^2\beta u^{\beta-1}|\nabla u|^2=\int f\eta^2u^{\beta+1}.
$$
By Young inequality and H\"older inequality:
\begin{align}\label{ie1}
|\beta|\int_{B_1} \eta^2 u^{\beta-1} |\nabla u|^2 \leq \frac{C}{|\beta|}\int_{B_1} |\nabla \eta|^2 u^{\beta+1}+(\int_{B_1} |f|^\frac{n}{2})^\frac{2}{n}(\int_{B_1} (\eta^2 u^{\beta+1})^\frac{n}{n-2})^\frac{n-2}{n}.
\end{align}
Applying the Sobolev inequality and Poincar\'e inequality  to $\eta u^\frac{\beta+1}{2}$
, we get
\begin{align*}
(\int_{B_1} (\eta u^\frac{\beta+1}{2})^\frac{2n}{n-2})^{\frac{n-2}{n}}&\leq \alpha_n\int_{B_1}|\nabla(\eta u^{\frac{\beta+1}{2}})|^2
\\&\leq 2\alpha_n \int_{B_1}(\nabla \eta)^2 u^{\beta+1}+2\alpha_n \int_{B_1}(\eta)^2|\nabla u^\frac{\beta+1}{2}|^2,
\end{align*}
which together with (\ref{ie1}) gives
\begin{align}\label{ie2}
\frac{4|\beta|}{(\beta+1)^2}\int_{B_1} \eta^2|\nabla u^\frac{\beta+1}{2}|^2 \leq (\frac{C}{|\beta|}+C\varepsilon_0)\int_{B_1} |\nabla\eta|^2 u^{\beta+1}+C\varepsilon_0\int_{B_1} \eta^2|\nabla u^\frac{\beta+1}{2}|.
\end{align}
When
$$
C\varepsilon_0\leq \frac{2|\beta|}{(\beta+1)^2},
$$
we have
$$
\frac{2|\beta|}{(\beta+1)^2}\int_{B_1} \eta^2|\nabla u^\frac{\beta+1}{2}|^2 \leq (\frac{C}{|\beta|}+\frac{2|\beta|}{(\beta+1)^2})\int_{B_1} |\nabla\eta|^2 u^{\beta+1},
$$
and
$$
\frac{2|\beta|}{(\beta+1)^2}\int_{B_1} |\nabla\eta u^\frac{\beta+1}{2}|^2 \leq (\frac{C}{|\beta|}+\frac{6|\beta|}{(\beta+1)^2})\int_{B_1} |\nabla\eta|^2 u^{\beta+1}.
$$
Take $\frac{1}{2}\leq r_1< r_2\leq 1$. Let $\eta\in C^\infty_0(B_{r_2})$, $\eta\equiv 1$ on $B_{r_1}$ and $|\nabla \eta|\leq\frac{C}{|r_2-r_1|}$.
By Poincar\'e inequality and Sobolev inequality, we get
$$
(\int_{B_{r_1}} |u^{\frac{\beta+1}{2}}|^{2^*})^\frac{1}{2^*}\leq C\left(\frac{(\beta+1)^2}{\beta^2}+1\right)\frac{1}{|r_2-r_1|}(\int_{B_{r_2}} (u^\frac{1+\beta}{2})^2)^\frac{1}{2},
$$
where $2^*=2\frac{n}{n-2}$.

Next we  deduce an uniform bound for $\|u\|_{L^p}$.

Let $\frac{\beta+1}{2}=\alpha$. We can choose $\varepsilon_0$ such that
$\|u\|_{L^{2^*\alpha}}<C$. Then by setting $\frac{\beta+1}{2}=2^*\alpha$ we can get $\|u\|_{L^{2^*\cdot 2^* \alpha}}<C$.  After several
iterations, we get an estimate of $\|u\|_{L^\frac{n}{n-2}}$.
So without loss of generality, we assume $\|u\|_{L^\frac{n}{n-2}}<C$.

Denote $\alpha=\frac{n}{n-2}$ and take
$$
\frac{n}{n-2}\geq p_0>1.
$$
Then
$$
(\int_{B_{r_1}} |u^{p_0\frac{\alpha(\beta+1)}{p_0}}|)^{\frac{p_0}{\alpha(\beta+1)}\frac{\beta+1}{2p_0}}\leq C\left(\frac{(\beta+1)^2}{\beta^2}+1\right)\frac{1}{|r_2-r_1|}(\int_{B_{r_2}} u^{p_0\frac{1+\beta}{p_0}})^{\frac{p_0}{\beta+1}\frac{\beta+1}{2p_0}},
$$
i.e.
\begin{equation}\label{moser1}
(\int_{B_{r_1}} |u^{p_0\frac{\alpha(\beta+1)}{p_0}}|)^{\frac{p_0}{\alpha(\beta+1)}}\leq \left(C(\frac{(\beta+1)^2}{\beta^2}+1)\frac{1}{|r_2-r_1|}\right)^\frac{2p_0}{\beta+1}(\int_{B_{r_2}} u^{p_0\frac{1+\beta}{p_0}})^{\frac{p_0}{\beta+1}}.
\end{equation}
Take $\beta+1=\alpha^mp_0$, $r_1=\frac{1}{2}+\frac{1}{2^{m+2}}$ and $r_2=\frac{1}{2}+\frac{1}{2^{m+1}}$, where $m=0,1,2,\cdots$, $m_0$, and
$$
m_0=\max\{m:C\varepsilon_0\leq \frac{2(\alpha^mp_0-1)}{(\alpha^mp_0)^2}\}.
$$
Rewrite \eqref{moser1} as follows:
$$
\|u^{p_0}\|_{L^{\alpha^{m+1}}(B_{r_1})}\leq C^\frac{2m}{\alpha^m}\|u^{p_0}\|_{L^{\alpha^m}(B_{r_2})},
$$
which implies
$$
\|u^{p_0}\|_{L^{\alpha^{m_0+1}(B_{\frac{1}{2}})}}\leq C^{\sum\limits_{i=0}^{+\infty} i\alpha^{-i}}\|u^{p_0}\|_{L^{1}(B_{1})}.
$$
Given $p\geq p_0$,  select $m_0$ such that $p<p_0\alpha^{m_0+1}$ and choose $\varepsilon_0$ under additional assumption:
$$
C\varepsilon_0\leq \min\{\frac{2(\alpha^mp_0-1)}{(\alpha^mp_0)^2}:m=0,1,\cdots,m_0\}.
$$
It follows
$$
\|u\|_{L^p(B_{\frac{1}{2}})}\leq C\|u\|_{L^{p_0\alpha^{m+1}}(B_{\frac{1}{2}})}\leq C\|u\|_{L^{p_0}(B_1)}\leq C.
$$

Now we return to the elliptic equation \eqref{equation.epsilon}.
For any $q<\frac{n}{2}$, we have
$$
(\int_{B_\frac{1}{2}} (fu)^q)^{\frac{1}{q}}\leq (\int_{B_\frac{1}{2}} f^\frac{n}{2})^{\frac{2}{n}}(\int_{B_\frac{1}{2}} u^\frac{n}{n-2q})^{\frac{n-2q}{n}}.
$$
Thus, if $p>\frac{n}{n-2q}$, we get
$$
\|u\|_{W^{2,q}(B_\frac{1}{4})}<C.
$$

Finally, we derive the estimate of $\|u^{-1}\|_{W^{2,q}}$.
Similar to the above arguments,  one can get $\|u^{-1}\|_{L^p(B_\frac{1}{4})}<C$. The estimate of $\|u^{-1}\|_{W^{2,q}}$ follows
from the following:
$$
\nabla u^{-1}=\frac{\nabla u}{u^2},\s \nabla^2u^{-1}=
\frac{\nabla^2u}{u^2}-2\frac{|\nabla u|^2}{u^3}.
$$
Since
$$
\nabla\log u=\frac{\nabla u}{u},\s \nabla^2\log u=
\frac{\nabla^2u}{u}-\frac{|\nabla u|^2}{u^2},
$$
we get the estimate of $\log u$.

Notice that for any positive constant $c$, $cu$ still satisfies the equation.
So we can get the estimate of $\|\log u\|_{W^{2,p}}$ without the assumption
that $\int_{B_1}\log u=0$.
\endproof

\subsection{Proof of Theorem \ref{main2}}The main goal of
this subsection is to prove Theorem \ref{main2}.
For any $x\in M$, on a small ball $B_r(x)\subset M$, $g_0|_{B_r(x)}$ can be regarded as a metric over $B_r\subset \R^n$ and we have the following equation
$$
-\frac{4(n-1)}{n-2}\Delta_{g_0} u_k=(-R(g_0)+R(g)u_k^\frac{4}{n-2})u_k.
$$
By Lemma \ref{Lalpha}, $\|\nabla\log u_k\|_{L^2(M)}<C$.  Let $-\log c_k$
be the mean value of $\log u_k$ over $M$.  By the Poincar\'e inequality,
$\|\log c_ku_k\|_{L^1(M)}<C$. Note $c_ku_k$ also satisfies
$$
-\frac{4(n-1)}{n-2}\Delta_{g_0}c_k u_k=(-R(g_0)+R(g)u_k^\frac{4}{n-2})c_ku_k.
$$
Cover $M$ by finitely many balls $B_{r_1}(x_1)$, $\cdots$,
$B_{r_m}(x_m)$. Assume for each $B_{r_i}(x_i)$,
$$
\int_{B_{2r_i}(x_i)} |R_k|^2dV_{g_k}<\varepsilon_0.
$$
Applying Lemma \ref{regularity} to $c_ku_k$, we get
 $\|c_ku_k\|_{W^{2,q}(M)}+\|(c_ku_k)^{-1}\|_{W^{2,q}(M)}<C$.

Since
$$
1=\int_Mu_k^\frac{2n}{n-2}=\frac{1}{c_k^\frac{2n}{n-2}}\int_M(c_ku_k)^\frac{2n}{n-2}<C\frac{1}{c_k^\frac{2n}{n-2}},
$$
and
\begin{eqnarray*}
\mathrm{Vol}^2(M,g_0)&\leq\int_Mu_k^\frac{2n}{n-2} dV_{g_0}\int_Mu_k^{-\frac{2n}{n-2}} dV_{g_0}=\int_Mu_k^{-\frac{2n}{n-2}} dV_{g_0}\\
&=c_k^\frac{2n}{n-2}\int_M(c_ku_k)^{-\frac{2n}{n-2}} dV_{g_0}
\leq Cc_k^\frac{2n}{n-2},
\end{eqnarray*}
we get the bound of $c_k$.
This proves the first part of Theorem \ref{main2}.

By Theorem \ref{main1}, we may assume the sequence $\{d_{g_k}\}$ converges to a distance function
$d$ with $d\leq d_g$. To finish the proof of Theorem \ref{main2},
we need to show $d\geq d_g$.
The key observation is the following:

\begin{lem}\label{du0/d0}
For any $\varepsilon$, we can find $\beta$ and $\tau$, which  only depend on $\varepsilon$, such that if
$$
\mu(B_{2\delta}(0))<\tau, \s \delta<\delta_0,
$$
then
$$
\frac{d_g(x,y)}{d(x,y)}\leq1+\varepsilon,\s \forall x, y\in B_{\beta\delta}^{g_0}(0).
$$
\end{lem}

\proof
We argue by contradiction and assume the lemma is not true. Then we can find $\delta_m\to0$ and $y_m, x_m\in B_{\delta_m}(0)$, such that
$\frac{|x_m-y_m|}{\delta_m}\rightarrow 0$ and
$$
\limsup_{k\rightarrow+\infty}\int_{B_{\delta(0)}}|R(g_k)|^\frac{n}{2}dV_{g_k}= 0,\s \frac{d(y_m,x_m)}{d_g(y_m,x_m)}\rightarrow a<1.
$$

For any fixed $m$, by the Sobolev Embedding Theorem and Lemma \ref{regularity}, and taking $q$ sufficiently closed to $\frac{n}{2}$, we can find $k_m$ such that $\{u_{k_m}\}$ and $\{\frac{1}{u_{k_m}}\}$ converges to $u$ and $\frac{1}{u}$ respectively in $W^{1,p_0}$, with $p_0<\frac{nq}{n-q}$. Moreover, using the H\"older's Inequality, we can get that, after passing to a subsequence, $\{\frac{u_{k_m}}{u}\}$ converges to $1$ in $W^{1,p}_{loc}(B_{\delta_0}(0))$, for some $p\in (n-1,p_0)$. So we have
$$
\left|\frac{d_{g_{k_m}}(y_m,x_m)}{d_g(y_m,x_m)}-\frac{d(y_m,x_m)}{d_g(y_m,x_{m})}\right|<\frac{1}{m},
$$
and
\begin{equation}\label{L2}
\int_{B_{\delta_m}(x_m)}\left|\nabla\frac{u_{k_m}(y)}{u(y)}\right|^pdy
+\int_{B_{\delta_m}(x_m)}\left|\frac{u_{k_m}(y)}{u(y)}-1\right|^pdy<\frac{1}{m},
\end{equation}
where $r_m=|y_m-x_m|$ and $B_{\delta_m}(x_m)\subset B_{\delta_0}(0)$. By (\ref{L2}), we have
\begin{equation}\label{L1}
r_m^{n-p}\int_{B_{\frac{\delta_m}{r_m}}(0)}\left|\nabla\frac{u_{k_m}(r_mx+x_m)}{u(r_mx+x_m)}\right|^pdx
+r_m^{n}\int_{B_{\frac{\delta_m}{r_m}}(0)}\left|\frac{u_{k_m}(r_mx+x_m)}{u(r_mx+x_m)}-1\right|^pdx<\frac{1}{m},
\end{equation}

For simplicity, we set $y_m=x_m+r_m(1,0,\cdots,0)$ in local coordinates. More precisely, we can take a segment $\gamma_m$ joining $x_m$ and $y_m$, such that
$$
\gamma_m(t)=x_m+r_m\gamma(t),\s where\s \gamma(t)=(t,0,\cdots,0),\s t\in[0,1].
$$
Let $u_m'=c_mu_{k_m}(x_m+r_mx)$, where $c_m$ is chosen such that
$$
0=\int_{B_\frac{1}{2}}\log u_m'.
$$
In the local coordinate, set $h_m(x)=(g_0)_{ij}(x_m+r_mx)dx^i\otimes dx^j$, which converges to $g_{\R^n}$ smoothly. Let $g_m'(x)=(u_m'(x))^\frac{4}{n-2}h_m(x)=\frac{c_m^\frac{4}{n-2}}{r_m^2}g_{k_m}(x_m+r_mx)$. Since for any $R>0$
$$
\int_{B_R(0)}|R(g_m')|^\frac{n}{2}dV_{g_m'}\leq\int_{B_{\delta_0(0)}}|R(g_k)|^\frac{n}{2}
dV_{g_k}\rightarrow 0.
$$
By Lemma \ref{regularity}, we may assume $u_m'$ converges to a positive harmonic function $u'$ weakly in $W^{2,q}_{loc}(\R^n)$
with  $\int_{B_\frac{1}{2}}\log u'=0$, for some $q<\frac{n}{2}$.
By Liouville's theorem, $u'$ is a constant. Since $\int_{B_\frac{1}{2}}\log u'=0$,
$u'=1$.
Hence, for a fixed $R>1$, $u_m'$ converges to $1$ in $W^{1,p}(B_R(0))$. By \eqref{L1}, $\frac{u(x_m+r_mx)}{u_{k_m}(x_m+r_mx)}$ converges to 1 in $W^{1,p}(B_R(0))$. According to the H\"older's Inequality, $u_m'(x)\frac{u(x_m+r_mx)}{u_{k_m}(x_m+r_{m}x)}$ converges to 1 in $W^{1,p'}(B_R(0))$, for some $p'\in (n-1,p)$. Then we have
\begin{align*}
\lim_{m\to \infty}c_m^\frac{2}{n-2}d_g(x_m,y_m)&\leq \lim_{m\to \infty}\int_{\gamma_m}(c_mu(\gamma_m(t)))^\frac{2}{n-2}
\\ &= \lim_{m\to \infty}\int_0^1(c_mu(\gamma_m(t)))^\frac{2}{n-2}
\sqrt{(g_0)_{ij}(\gamma_m(t))\dot{\gamma}_m^i(t)\dot{\gamma}_m^j(t)}dt
\\ &=\lim_{m\to \infty} \int_0^1 \left(c_mu_{k_m}(x_m+r_{m}\gamma(t))\frac{u(x_m+r_m\gamma(t))}{u_{k_m}(x_m+r_{m}\gamma(t))}\right)^\frac{2}{n-2}dt
\\ &=\lim_{m\to \infty} \int_0^1\left(u_m'(\gamma(t))\frac{u(x_m+r_m\gamma(t))}{u_{k_m}(x_m+r_{m}\gamma(t))}\right)^\frac{2}{n-2}dt.
\end{align*}
Therefore, we get
$$
\lim_{m\rightarrow+\infty} c_m^\frac{2}{n-2}d_g(x_m,y_m)\leq 1.
$$
By Proposition \ref{complete.case},
$$
c_m^\frac{2}{n-2}d_{g_{k_m}}(x_m,y_m)=d_{g_m'}(0,(1,0\cdots,0))\rightarrow 1.
$$
Hence,
$$
\frac{d(x_m,y_m)}{d_g(x_m,y_m)}\geq1
$$
when $m$ is sufficiently large.

\endproof

{\it The proof of Theorem \ref{main2}:} It remains to show that $d_g\leq d$. Let $\varepsilon$, $\tau$ and $\beta$ be in Lemma \ref{du0/d0} and set $A_\tau=\{x:\mu(\{x\})>\tau\}$. Obviously, $A_\tau$
is a finite set.   Then, for any $\delta>0$,  we have
$$
\int_{B_\delta(x)}|R(g_k)|^\frac{n}{2}dV_{g_k}<\tau,\s whenever\s B_{2\delta}(x)
\cap A_\tau=\emptyset,
$$
when $k$ is sufficiently large. It follows that
$$
\frac{d_g(x,y)}{d(x,y)}<1+\varepsilon,
$$
whenever $d_{g_0}(x,y)<\beta\delta$ and $x\notin B_\delta(A_\tau)$.

We say a metric space is a length space if, for each $x$
and $y$ in this space, there exists a minimal geodesic joining them (see \cite[p.148]{Fukaya}).
Since $(M,d)$ is also the Gromov-Haudorff limit of $(M,d_{g_k})$, $(M,d)$ is a length space according to \cite [Proposition 1.10]{Fukaya}.  Let $\gamma$ be the minimal geodesic defined in $(M,d)$ joining
$x_1$ and $x_2$, i.e. $\gamma:[0,a]\rightarrow (M,d)$
is a continuous map which satisfies
$$
\gamma(0)=x_1,\s \gamma(a)=x_2,\s and\s
d(\gamma(s),\gamma(s'))=|s-s'|,\s\forall s,s'\in[0,a].
$$
We claim that $\gamma$ is also continuous in $(M,g_0)$. For otherwise, we can find $t_k\rightarrow t$ and $a>0$, such that
$d_{g_0}(\gamma(t_k),\gamma(t))>a$.
By \eqref{positive.distance2}, there exists $a'>0$, such that
$$
|t_k-t|\geq d(\gamma(t_k),\gamma(t))>a',
$$
which is impossible.

Now we consider two cases. The first case is when $\gamma\cap A_\tau=\emptyset$, say $d(A_\tau,\gamma[0,a])>0$.
Since $\gamma$ is continuous,  we may assume
$$
d_{g_0}(A_\tau,\gamma[0,a])>\delta>0.
$$
Then there exists
$$
s_0=0<s_1<\cdots<s_m=a,
$$
such that
$$
d_{g_0}(\gamma(s_{i+1}),\gamma(s_i))<\beta\delta.
$$
It follows that
\begin{eqnarray}\label{d0.d0'}\nonumber
d(x_1,x_2)&=&\sum_{i=0}^{m-1} d(\gamma(s_i),\gamma(s_{i+1}))\\
&\geq& (1+\varepsilon)^{-1}\sum_{i=0}^{m-1}d_g(\gamma(s_i),
\gamma(s_{i+1}))\\\nonumber
&\geq& (1+\varepsilon)^{-1}d_g(x_1,x_2)
\end{eqnarray}

The remaining case is when $\gamma\cap A_\tau\neq
\emptyset$.  Let
$$
\gamma\cap A_\tau=\{\gamma(a_1), \cdots, \gamma(a_i)\}.
$$
The distance function $d$ is bounded by
\begin{eqnarray*}
d(x_1,x_2)&\geq& d(x_1,\gamma(a_1-\varepsilon'))
+d(\gamma(a_1+\varepsilon'),\gamma(a_2-\varepsilon'))+
\cdots+d(\gamma(a_i+\varepsilon',x_2))\\
&\geq&
(1+\varepsilon)^{-1}(d_g(x_1,\gamma(a_1-\varepsilon'))
+
\cdots+d_g(\gamma(a_i+\varepsilon',x_2))).
\end{eqnarray*}
When $\varepsilon'\rightarrow 0$, we get \eqref{d0.d0'} again.

Finally, by letting $\varepsilon\rightarrow 0$, we get the desired inequality.
\endproof
{\small

\end{document}